\documentclass[12pt]{article}
\usepackage{amssymb,latexsym,amscd,amsmath,amsfonts,enumerate,supertabular}
\usepackage{graphicx}
\usepackage{tabularx}
\usepackage{caption}
\usepackage{epstopdf}
\usepackage{color}
\usepackage{graphicx}

\newtheorem{Theorem}{Theorem}[section] 
\newtheorem{Proposition}[Theorem]{Proposition}

\begin{document}
\title{{\bf High order numerical methods  for solving high orders functional differential equations }}

\author{ Dang Quang A$^{\text {a,b}}$,  Dang Quang Long$^{\text c}$\\
$^{\text a}$ { Center for Informatics and Computing, VAST}\\
{ 18 Hoang Quoc Viet, Cau Giay, Hanoi, Vietnam}\\
{ Email: dangquanga@cic.vast.vn}\\
$^{\text b}$ { Thai Nguyen University of Information and Communication }\\ 
{Technology, Thai Nguyen, Viet Nam}\\
$^{\text c}$ { Institute of Information Technology, VAST,}\\
{ 18 Hoang Quoc Viet, Cau Giay, Hanoi, Vietnam}\\
{ Email: dqlong88@gmail.com}}
\date{ }

\date{}          
\maketitle
\begin{abstract}
In this paper we construct high order numerical methods for solving third and fourth orders nonlinear functional differential equations (FDE). They are based on the discretization of  iterative methods on continuous level with the use of the trapezoidal quadrature formulas with corrections. Depending on the number of terms in the corrections we obtain methods of $O(h^4)$ and $O(h^6)$ accuracy. Some numerical experiments demonstrate the validity of the obtained theoretical results. 
   The approach used here for the third and fourth orders nonlinear functional differential equations
can be applied to functional differential equations of any orders.

\end{abstract}
\noindent {\bf Keywords: }High order functional differential equations;  Iterative method;  High order numerical method.
%

\section{Introduction}\label{sec1}
Recently, in \cite{Dang1} we established the existence and uniqueness results and constructed a numerical method of second order of accuracy for solving the FDE
\begin{equation}\label{eq1}
u^{\prime \prime \prime}=f(t,u(t), u(\varphi (t))), \quad t\in (0,a)
\end{equation}
associated with the general linear-type two-point boundary conditions. The method is based on the use of Green function and trapezoidal formula for computing integrals at each iteration. We proved that the method is of $O(h^2)$ accuracy. In the last year, Bica and Curila \cite{Bica1} constructed successive approximations for equivalent integral equation with the use of cubic spline interpolation at each iterative step for solving the equation \eqref{eq1} in $[0, 1]$ with the boundary conditions
\begin{equation}\label{eq2}
u(0)=c_1, u'(0)=c_2, u'(1)=c_3.
\end{equation}
They prove that the maximal order of convergence of the method is $O(h^3)$.  The proof of this result is very complicated. It should be said that the method used in \cite{Bica1} is the same method of iterated splines proposed in \cite{Bica2} for even order FDEs. There are also numerical examples which show that the method is of order $O(h^2)$ for second order FDE and of order $O(h^4)$ for fourth order FDE. In \cite{Bica2,Bica3} the authors proved the maximal order of convergence of the method for the fourth order FDE with the clamped boundary conditions is $O(h^4)$. Also, in \cite{Bica4} the authors used Catmull-Rom cubic spline interpolation at each iteration of the Green function technique to third order FDE and come to the same conclusion that the maximal order of convergence is $O(h^3)$.\par 
It should be emphasized that the approach to boundary value problems (BVPs) in \cite{Bica1,Bica2,Bica4} is based on the discretization of the iterative method applied to the equivalent nonlinear integral equation
\begin{equation*}
u(t)=g(t)+\int_0^a G(t,s) f(s,u(s),u(\varphi(s)))ds.
\end{equation*}
Due to this approach, with the use of the cubic spline interpolation at each iteration the construction of the numerical algorithms and the proof of their convergence are very cumbersome and complicated.\par 
Differently from Bica, in \cite{Dang1} we constructed numerical method for solving the third order FDE \eqref{eq1} on the base of discretization of the iterative method on continuous level resulting from the operator equation for nonlinear term (or right hand side) of the equation. We proved that the method obtained by using trapezoidal quadrature formula at each iteration is of $O(h^2)$.\par 
In this paper, using trapezoidal quadrature formula with corrections we construct numerical methods of order $O(h^4)$ and $O(h^6)$ for the equation \eqref{eq1} with the boundary conditions \eqref{eq2}. We also construct numerical methods of high orders for a BVP for fourth order FDE. The technique  used here is an further development of our method 
 for designing high orders numerical methods for solving high orders ODEs \cite{Dang4,Dang5}.
Some numerical examples confirm the theoretical results.

\section{Preliminaries}
\subsection{Iterative method on continuous level for FDE}
In this subsection we recall the iterative method on continuous level for solving the equation 
\begin{equation}\label{eq3}
u^{\prime \prime \prime}=f(t,u(t), u(\varphi (t))), \quad t\in (0,a)
\end{equation}
with the general linear boundary conditions
\begin{equation}\label{eq4}
 \begin{split}
 B_1[u]=\alpha_1 u(0)+\beta_1 u'(0) + \gamma_1 u''(0) =b_1,\\
 B_2[u]=\alpha_2 u(0)+\beta_2 u'(0) + \gamma_2 u''(0) =b_2,\\
 B_3[u]=\alpha_3 u(1)+\beta_3 u'(1) + \gamma_3 u''(1) =b_3,\\
 \end{split}
 \end{equation}
or
 \begin{equation}\label{eq5}
 \begin{split}
 B_1[u]=\alpha_1 u(0)+\beta_1 u'(0) + \gamma_1 u''(0) =b_1,\\
 B_2[u]=\alpha_2 u(1)+\beta_2 u'(1) + \gamma_2 u''(1) =b_2,\\
 B_3[u]=\alpha_3 u(1)+\beta_3 u'(1) + \gamma_3 u''(1) =b_3,\\
 \end{split}
 \end{equation}
%

The results of existence and uniqueness of the problem \eqref{eq3}-\eqref{eq4} are given in \cite{Dang1}. For solving the problem we proposed the following iterative method:
\begin{enumerate}
\item Given 
\begin{equation}\label{iter1c}
\psi_0(t)=f(t,0,0).
\end{equation}
\item Knowing $\psi_k(t)$  $(k=0,1,...)$ compute
\begin{equation}\label{iter2c}
\begin{split}
u_k(t) &= g(t)+\int_0^a G(t,s)\psi_k(s)ds ,\\
v_k(t) &= g(\varphi (t))+\int_0^a G(\varphi(t),s)\psi_k(s)ds .
\end{split}
\end{equation}
\item Update
\begin{equation}\label{iter3c}
\psi_{k+1}(t) = f(t,u_k(t),v_k(t)),
\end{equation}
\end{enumerate}
where $g(t)$ is the polynomial of second degree satisfying the boundary conditions, $G(t,s)$ is the Green function of the problem. To state the theorem of convergence of the above iterative method we recall some notations:\\
\begin{equation}\label{eqM0}
M_0  =\max_{0\leq t\leq a} \int_0^a |G(t,s)|ds ,
\end{equation}
\begin{equation}\label{eqDM}
\mathcal{D}_M= \Big \{ (t,u,v) \mid  0\leq t\leq a; |u|\leq \|g\|+M_0 M; |v|\leq \|g\|+M_0 M \Big \},
\end{equation}
where $M$ is any  positive number, $\|g\|= \max _{0\le t\le a}|g(t)|$.\\
As usual, we denote by $B[0,M]$ the closed ball of the radius $M$ centered at $0$ in the space of continuous functions $C[0,a]$. 

\begin{Theorem}\label{theorem1} (see Theorems 2.2 and 3.1 in \cite{Dang1})
Assume that:
\begin{description}
\item (i) The function $\varphi (t)$ is a continuous map from $[0,a]$ to $[0,a]$.
\item (ii) The function $f(t,u,v)$ is continuous and bounded by $M$ in the domain $\mathcal{D}_M$, i.e.,
\begin{equation}\label{eq:12}
|f(t,u,v)|\le M \quad \forall (t,u,v) \in \mathcal{D}_M.
\end{equation}
\item (iii) The function $f(t,u,v)$ satisfies the Lipschitz conditions in the variables $u,v$ with the coefficients $L_1, L_2 \ge 0$ in $\mathcal{D}_M$, i.e.,
\begin{equation}\label{eq:13}
\begin{split}
|f(t,u_2,v_2)-f(t,u_1,v_1)|\le L_1 |u_2-u_1|+L_2 |v_2-v_1| \quad \\
\forall (t,u_i,v_i)\in \mathcal{D}_M\;
(i=1,2)
\end{split}
\end{equation}
\item (iv)   
\begin{equation}\label{eq:q}
q:= (L_1+L_2)M_0 <1.
\end{equation}
\end{description}
Then:\\
i) The the problem \eqref{eq3}-\eqref{eq4} has a unique solution $u(t) \in C^3[0,a]$, satisfying the estimate
\begin{equation}\label{eq:14}
|u(t)| \le \|g\| +M_0M \quad \forall t\in[0,a].
\end{equation}
ii) The above iterative method converges and there holds the estimate
\begin{equation*}
\|u_k-u\| \leq M_0p_kd, 
\end{equation*} 
where $u$ is the exact solution of the problem and 
\begin{equation}\label{eqpkd}
p_k=\dfrac{q^k}{1-q} ,\; d=\|\psi _1 -\psi _0\|
\end{equation}
\end{Theorem}
\subsection{Trapezoidal quadrature formula with corrections}
Let $h=a/n,$ where $N$ is a some positive integer, and $s_i=(i-1)h, i=1,\dots , n+1$. Then the Euler-Maclaurin formula has the form (see \cite{Hild})
\begin{equation}\label{eqHild}
\int _0^a \Phi (s) ds = T_{\Phi}(h) -\sum_{l=1}^{p-1} \frac{B_{2l}}{(2l!)} \Big ( \Phi ^{(2l-1)}(a)-\Phi ^{(2l-1)}(0) \Big ) +O(h^{2p}),
\end{equation}
where $\Phi \in C^{2p}[0, a]$, $B_{2l}$ are Bernoulli numbers,
\begin{equation}
T_{\Phi}(h)=\frac{h}{2}(\Phi _1 +\Phi _{n+1})+ \sum _{i=2}^{n} h \Phi_{i}, 
\end{equation}
with $ \Phi_i=\Phi (s_i).$
Now, for a fixed $t\in [0, a]$ let
\begin{equation*}
\Phi (s)= g(t,s) \psi (s),
\end{equation*}
where $g(t,s)$ is continuous in $[0, a] \times [0, a]$, and may have discontinuous derivatives at the point $s=t$, meanwhile $\psi (s) \in C^{2p}[0, a]$. Then using the above Euler-Maclaurin formula Sidi and Pennline \cite{Sidi} obtained the following formula
\begin{equation}\label{Euler-Mac}
\begin{split}
\int _0^a \Phi (s) ds =  &T_{\Phi}(h) -\sum_{l=1}^{p-1} \frac{B_{2l}}{(2l!)} \Big \{ [ \Phi ^{(2l-1)}(a)-\Phi ^{(2l-1)}(0) ] \\
 &-[ \Phi ^{(2l-1)}(t^+)-\Phi ^{(2l-1)}(t^-)  ]\Big \} +O(h^{2p}),
\end{split}
\end{equation}
where $\psi (t^{+})$ and $\psi (t^{-})$ are the one-sided limits of the function $\psi (s)$ at $t$.
In the particular case $p=2$ we have
\begin{equation}\label{eqCasep2}
\int _0^a \Phi (s) ds =  T_{\Phi}(h) -\frac{h^2}{12} \Big \{ [ \Phi '(a)-\Phi '(0) ]
  -[ \Phi '(t^+)-\Phi '(t^-)  ]\Big \} +O(h^4)
\end{equation}
and for the case $p=3$
\begin{equation}\label{eqCasep3}
\begin{split}
\int _0^a \Phi (s) ds =  T_{\Phi}(h) -\frac{h^2}{12} \Big \{ [ \Phi '(a)-\Phi '(0) ]
  -[ \Phi '(t^+)-\Phi '(t^-)  ]\Big \} \\
  +\frac{h^4}{720} \Big \{ [ \Phi '''(a)-\Phi '''(0) ]
  -[ \Phi '''(t^+)-\Phi '''(t^-)  ]\Big \} + O(h^6).
\end{split}
\end{equation}
If the function $\Phi (s)$ has a jump at  point $t\in (0,a)$ then the above formulas instead $T_{\Phi}(h)$ it should be $T_{\Phi^{*}}(h)$, where 
\begin{equation*}
\Phi^{*}(s)=\left\{\begin{array}{ll}
\Phi(s), \quad s\neq t,\\
\dfrac{1}{2}[\Phi(t^{+})+ \Phi(t^{-}) ], \quad s=t
\end{array}\right.
\end{equation*}

\section{Construction of numerical methods of orders 4 and 6 for the third order FDE}
\subsection{Computing the integrals of the type in \eqref{iter2c}}
Now we apply the formulas \eqref{eqCasep2} and \eqref{eqCasep3} to construct numerical methods of orders 4 and 6. Recall that the Green function associated with the problem \eqref{eq1}  with respect to $t$ have the form
\begin{equation}\label{eqG0}
\begin{aligned}
G(t,s)=\left\{\begin{array}{ll}
\dfrac{s}{2}( \dfrac{t^2}{a}-2t+s), \quad 0\le s \le t \le a,\\
\, \, \dfrac{t^2}{2}(\dfrac{s}{a}-1), \quad 0\le t \le s \le a.\\
\end{array}\right.
\end{aligned}
\end{equation}
We have $G(t,0)=G(t,a)=0$ and
\begin{equation*}
\frac{\partial G(t,s) }{\partial s}=\left\{\begin{array}{ll}
\dfrac{t^2}{2a}-t+s), \quad 0\le s \le t \le a,\\
\dfrac{t^2}{2a}, \quad 0\le t \le s \le a,\\
\end{array}\right.
\end{equation*}
For a fixed $t$ put $$\Phi (s)= G(t,s) \psi (s),$$
where $\psi \in C^1[0, a]$.
It is easy to verify that
\begin{align*}
 &\Phi '(a)- \Phi '(0)= G_1(t,a)\psi(a) - G_1(t,0)\psi(0),\\
  &\Phi '(t^+)-\Phi '(t^-) =0,\\
  &\Phi '''(a)- \Phi '''(0)=\frac{3t^2}{2a}\psi '' (a) -3\psi '(0) -3 (\frac{t^2}{2a}-t) \psi ''(0),\\
  &\Phi '''(t^+)-\Phi '''(t^-) =-3 \psi '(t),
 \end{align*}
 where $G_1(t,s)=\frac{\partial G(t,s) }{\partial s}$. 
Therefore, from \eqref{eqCasep2} we have
\begin{equation}\label{Gpsi}
\begin{split}
&\int_0^a G(t,s)\psi(s)ds = T_{\Phi}(h) - \frac{h^2}{12} \Big \{ G_1(t,a)\psi(a) - G_1(t,0)\psi(0) \Big \}+ O(h^4),\\
&\int_0^a G(t,s)\psi(s)ds = T_{\Phi}(h) - \frac{h^2}{12} \Big \{ G_1(t,a)\psi(a) - G_1(t,0)\psi(0) \Big \}\\
&+ \frac{h^4}{720}  \Big \{\frac{3t^2}{2a}\psi '' (a) -3\psi '(0) -3 (\frac{t^2}{2a}-t) \psi ''(0) + 3 \psi '(t)  \Big \}+ O(h^6).
\end{split}
\end{equation}
At the points $t_i =(i-1)h, i=2,...,n+1$ we obtain
\begin{equation}
\int_0^1 G(t_i,s)\psi (s) ds = L_4(G, t_i) \psi + O(h^4),
\end{equation}
and
\begin{equation}
\int_0^1 G(t_i,s)\psi (s) ds = L_6(G, t_i) \psi + O(h^6),
\end{equation}
where 
\begin{equation}\label{L4Gti}
L_4(G, t_i) \psi = \sum_{j=1}^{n+1}h\rho_j G(t_i,t_j)\psi_j - \frac{h^2}{12} \Big \{ G_1(t_i,t_{n+1})\psi_{n+1}- G_1(t_i,t_{1})\psi_{1}\Big \}, \; i=2,...,n+1,
\end{equation}
\begin{equation}\label{L6Gti}
\begin{split}
&L_6(G, t_i) \psi = \sum_{j=1}^{n+1}h\rho_j G(t_i,t_j)\psi_j - \frac{h^2}{12} \Big \{ G_1(t_i,t_{n+1})\psi_{n+1}- G_1(t_i,t_{1})\psi_{1}\Big \}\\ 
& +\frac{h^4}{720}  \Big \{\frac{3t_i^2}{2a}D^{(2)}_2 \psi_{n+1} -3D^{(1)}_2\psi_1 
-3 (\frac{t_i^2}{2a}-t_i) D^{(2)}_2 \psi_1 + 3 D^{(1)}_2 \psi _i  \Big \} \; i=2,...,n+1,
\end{split}
\end{equation}
In the above formula  $\rho_j$ are the weights 
\begin{equation*}
\rho_j = 
\begin{cases}
1/2,\; j=0,N\\
1, \; j=1,2,...,N-1
\end{cases}
\end{equation*}
and for short we denote $\psi_j =\psi (t_{j}), j=1,...,n+1$. Besides, $D^{(m)}_k \psi_i$ denotes the difference approximation of $m-$derivative of the function $\psi (t)$ at point $t_i$ with the accuracy $O(h^k)$ (see \cite{Li}).\\
Since $\int_0^1 G(t_i,s)\varphi(s)ds =0$ at $t_1=0$ due to $G(0,s)=0$ we set 
\begin{equation}\label{L4Gt1}
L_4(G, t_1) \psi = 0,\; L_6(G, t_1) \psi = 0
\end{equation}

\subsection{Discrete iterative method for the BVP \eqref{eq3}-\eqref{eq4}}
Cover the interval $[0, a]$   by the uniform grid $\bar{\omega}_h=\{t_i=(i-1)h, \; h=1/n, i=1,2,...,n+1  \}$ and denote by $\Psi_k(t), U_k(t), V_k(t)$ the grid functions, which are defined on the grid $\bar{\omega}_h$ and approximate the functions $\psi_k (t), u_k(t), v_k(t)$ on this grid, respectively.\par
Consider the following iterative methods\\
{\bf Iterative Method 1 (Fourth order method): }\\
\begin{enumerate}[(i)]
\item Given 
\begin{equation}\label{iter1d}
\Psi_0(t_i)=f(t_i,0,0),\ i=1,...,n+1. 
\end{equation}
\item Knowing $\Psi_k(t_i),\; k=0,1,...; \; i=1,...,n+1, $  compute approximately the definite integrals \eqref{iter2c} by the trapezoidal formulas with corrections:
\begin{equation}\label{iter2d}
\begin{split}
U_k(t_i) &= g(t_i) + L_4(G,t_i)  \Psi_k, \\
V_k(t_i) &= g(\xi_i) + L_4(G,\xi_i)  \Psi_k,   i=1,...,n+1,
\end{split}
\end{equation}
where $\xi_i = \varphi(t_i)$ and $L_4(G,t_i)  \Psi_k$ is defined by \eqref{L4Gti}-\eqref{L4Gt1}by replacing the function $\psi$ on the grid by the grid function $\Psi_k$
\item Update
\begin{equation}\label{iter3d}
\Psi_{k+1}(t_i) = f(t_i,U_k(t_i),V_k(t_i)),\;  i=1,...,n+1.
\end{equation}
\end{enumerate}

\noindent {\bf Iterative Method 2 (Sixth order method): }\\
Notice that according to \eqref{Gpsi}
\begin{equation}
\begin{split}
&\int_0^a G(\xi_i,s)\psi(s)ds = T_{\Phi}(h) - \frac{h^2}{12} \Big \{ G_1(\xi_i,a)\psi(a) - 
G_1(\xi_i,0)\psi(0) \Big \}\\
&+ \frac{h^4}{720}  \Big \{\frac{3\xi_i^2}{2a}\psi '' (a) -3\psi '(0) -3 (\frac{\xi_i^2}{2a}-\xi_i) \psi ''(0) + 3 \psi '(\xi_i)  \Big \}+ O(h^6).
\end{split}
\end{equation}
We will approximate $\psi '(\xi_i)$ by a difference formula with the accuracy $O(h^2)$.\\
If $\xi_i =0$ we have at once $\int_0^a G(\xi_i,s)\psi(s)ds =0$ since $G(0,s)=0$. This occurs when $i=1$ for the functions  $\varphi (t) =\alpha t $ or $\varphi (t) =\alpha t^2, \; \alpha >0 $. Consider the case $\xi_i >0, \; (i=2,..., n+1)$. Suppose that there exists $m \ge 1$ such that $t_m \le \xi_i <t_{m+1}$. It is easy to see that 
$m=[\frac{\xi_i}{h}]+1$, where $[a]$ denotes the integer part of the number $a$.\\
By the Taylor expansion we have 
\begin{equation}\label{derPsi}
\psi '(\xi_i) = \psi '(t_m)+\psi ''(t_m)(\xi_i -t_m)+O(h^2).
\end{equation}
Approximating $\psi '(t_m)$ by $D^{(1)}_2\psi_m$ and $\psi ''(t_m)$ by  $D^{(2)}_2\psi_m$, i.e., 
\begin{align*}
\psi '(t_m) = D^{(1)}_2\psi_m +O(h^2), \; \psi ''(t_m) = D^{(2)}_2\psi_m +O(h^2)
\end{align*}
from \eqref{derPsi} we obtain
\begin{equation}
\psi '(\xi_i) = D^{(1)}_2\psi_m + D^{(2)}_2\psi_m (\xi_i -t_m)+ O(h^2).
\end{equation}
Now we define
\begin{equation}
\begin{split}
&L_6(G, \xi_i) \psi = \sum_{j=1}^{n+1}h\rho_j G(\xi_i,t_j)\psi_j - \frac{h^2}{12} \Big \{ G_1(\xi_i,t_{n+1})\psi_{n+1}- G_1(\xi_i,t_{1})\psi_{1}\Big \}\\ 
& +\frac{h^4}{720}  \Big \{\frac{3\xi_i^2}{2a}D^{(2)}_2 \psi_{n+1} -3D^{(1)}_2\psi_1 
-3 (\frac{\xi_i^2}{2a}-\xi_i) D^{(2)}_2 \psi_1\\
& + 3 (D^{(1)}_2\psi_l + D^{(2)}_2\psi_l (\xi_i -t_l)) \Big \}, \; i=2,...,n+1.
\end{split}
\end{equation}
Then we have 
\begin{equation}
\int_0^a G(\xi_i,s)\psi(s)ds = L_6(G, \xi_i) \psi +O(h^6).
\end{equation}
Now we are able to describe a sixth order method or Method 2. \\
Method 2 is the same as Method 1 with a difference that in the step (ii) instead of \eqref{iter2d} there are
\begin{equation}\label{iter2d2}
\begin{split}
U_k(t_i) &= g(t_i) + L_6(G,t_i)  \Psi_k, \\
V_k(t_i) &= g(\xi_i) + L_6(G,\xi_i)  \Psi_k,   i=1,...,n+1,
\end{split}
\end{equation}
Now we study the convergence of the above iterative methods.
\begin{Proposition}\label{prop4}
Under the assumptions of Theorem \ref{theorem1} and  additional assumption that  the function $f(t,x,y)$ is sufficiently smooth, for Iterative method $r \; (r=1, 2)$ we have 
\begin{equation}\label{eq:prop4a}
\|\Psi_k -\psi_k  \|= O(h^{2(r+1)}),\ \|U_k -u_k  \|=O(h^{2(r+1)}), \ k=0,1,...
\end{equation}
\noindent where $\|.\|_{C(\bar{\omega}_h)}$ is the max-norm of function on the grid $\bar{\omega}_h$.
\end{Proposition}
{\bf Proof.}  The proposition can be proved by induction in a similar way as Proposition 3.5 in \cite{Dang1} if taking into account the orders 4 and 6 of quadrature formulas used in design of the discrete methods and the linearity of $L_4(G,t_i)$ and $L_6(G,t_i)$ as an operators acting on grid function $\Psi_k$.\\
Combining the above proposition and  Theorem \ref{theorem1} we obtain the following theorem.
\begin{Theorem}\label{theorem2}
Under the assumptions of the above proposition and  Theorem \ref{theorem1} for the approximate solution of the problem \eqref{eq3}-\eqref{eq4} obtained by the discrete Iterative Method $r\; (r=1, 2)$  on the uniform grid with gridsize $h$ we have the estimates
\begin{equation*}
\begin{split}
\|U_k-u\| &\leq M_0p_kd +O(h^{2(r+1)}),  
\end{split}
\end{equation*}
where $M_0$ are defined by \eqref{eqM0} and $p_k, d$ are defined by \eqref{eqpkd}.
\end{Theorem}
\section{Examples for third order FDEs}
We test the proposed method on the examples taken from \cite{Bica1} for comparison and on some other examples for demonstrating the efficiency of the methods. The iterative methods are performed  until $ \|\Psi_{k+1}-\Psi_k \| \le 10^{-16}.$\\

\noindent {\bf Example 4.1.} (Example 1 in \cite{Bica1})
Consider the third order boundary value problem
\begin{align*}
u''' (t)=e^{-t}[u(t)]^{\frac{3}{2}} . u(\frac{t}{2}), t\in [0, 1]\\
u(0)=u'(0)=1,\ u'(1)=e
\end{align*}
with the exact solution $u(t)=e^{t}$.\\
The results of convergence of  Iterative method 1  are given in Table \ref{Tab1}, where $n+1$ is the number of grid points, $K$ is the number of iterations performed, $Order$ is the order of convergence calculated by the formula
$$ Order=\log _2 \frac{\|U^{n/2}_K-u\|}{\|U^{n}_K-u\|}.
$$
The superscripts $n/2$ and $n$ of $U_K$ mean that $U_K$ is computed on the grid with the corresponding number of grid points.
\begin{table}[!ht]
    \centering
     \caption{The convergence of Method 1 for Example 4.1 }
     \label{Tab1}
    \begin{tabular}{ccccc}
        \hline
       $ n$ & $K$ & $Error$ & $h^4$ & $Order$ \\ \hline
        8 & 16 & 3.5273e-06 & 2.4414e-04 & ~ \\ 
        16 & 16 & 2.2041e-07 & 1.5259e-05 & 4.0003 \\ 
        32 & 16 & 1.3775e-08 & 9.5367e-07 & 4.0001 \\ 
        64 & 17 & 8.6093e-10 & 5.9605e-08 & 4.0000 \\ 
        128 & 17 & 5.3808e-11 & 3.7253e-09 & 4.0000 \\ 
        256 & 17 & 3.3631e-12 & 2.3283e-10 & 4.0000 \\ 
        512 & 17 & 2.1050e-13 & 1.4552e-11 & 3.9979 \\ \hline
    \end{tabular}
\end{table} 
From Table \ref{Tab1} it is obvious that the iterative method has the accuracy order $O(h^4)$ and the convergence order is $4$. For comparison of the accuracy of the method with the results of \cite{Bica1} we perform computation for $n=10, 100, 1000$. The errors of approximate solution are 1.4447e-06, 1.4444e-10, 1.4655e-14, and the numbers of iterations are $16, 16, 17$, respectively. So, these results are almost the same in \cite{Bica1} although there it is proved theoretically that the convergence is only $O(h^3)$ .\\
The convergence of  Method 2 is reported in Table \ref{Tab1a}, from which it is seen that the accuracy of the method is $O(h^6).$ 

\begin{table}[!ht]
    \centering
    \caption{The convergence of Method 2 for Example 4.1 }
     \label{Tab1a}
    \begin{tabular}{cccccc}
    \hline
        $N$ & $h^6$ & Iter & Error & Iter of Bica & Error of Bica \\ \hline
        10 & 1.0000e-06 & 16 & 6.1215e-09 & 16 & 1.376e-6 \\ 
        20 & 1.5625e-08 & 16 & 1.5168e-10 & ~ & ~ \\ 
        40 & 2.4414e-10 & 18 & 1.2227e-11 & ~ & ~ \\ 
        80 & 3.8147e-12 & 17 & 8.1268e-13 & ~ & ~ \\ 
        100 & 1.0000e-12 & 16 & 3.3529e-13 & 16 & 1.377e-10 \\ 
        200 & 1.5625e-14 & 17 & 2.1760e-14 & ~ & ~ \\ \hline
    \end{tabular}
\end{table}

\noindent {\bf Example 4.2.} (Example 2 in \cite{Bica1})
Consider the third order boundary value problem
\begin{align*}
u''' (t)&=-\frac{4}{(t+1)^4}-([u(t)]^4+[u(t)]^3). u(\frac{t}{2}), t\in [0, 1]\\
u(0)&=1,\ u'(0)=-1, u'(1)=-\frac{1}{4}
\end{align*}
with the exact solution $u(t)=\dfrac{1}{t+1}$. \\
The results of computation by Method 1 for Example 2 are presented in Table \ref{Tab2}. 
For $n=10, 100, 1000$ the errors of approximate solution are 3.9977e-05, 4.0704e-09, 4.0723e-13, and the numbers of iterations are $16, 17, 18$, respectively.\\
\begin{table}[!ht]
    \centering
     \caption{The convergence of Method 1 for Example 4.2 }
     \label{Tab2}
    \begin{tabular}{ccccc}
        \hline
       $ n$ & $K$ & $Error$ & $h^4$ & $Order$ \\ \hline
        8 & 16 & 9.6622e-05 & 2.4414e-04 & ~ \\ 
        16 & 16 & 6.1678e-06 & 1.5259e-05 & 3.9695 \\ 
        32 & 16 & 3.8756e-07 & 9.5367e-07 & 3.9923 \\ 
        64 & 18 & 2.4255e-08 & 5.9605e-08 & 3.9981 \\ 
        128 & 17 & 1.5165e-09 & 3.7253e-09 & 3.9995 \\ 
        256 & 17 & 9.4787e-11 & 2.3283e-10 & 3.9999 \\ 
        512 & 18 & 5.9244e-12 & 1.4552e-11 & 3.9999 \\ \hline
    \end{tabular}
\end{table}

The convergence of  Method 2 for Example 4.2 is reported in Table \ref{Tab2a}, from which it is seen that the accuracy of the method is $O(h^6).$ 

\begin{table}[!ht]
    \centering
        \caption{The convergence of Method 2 for Example 4.2 }
     \label{Tab2a}
    \begin{tabular}{cccccc}
    \hline
 $N$ & $h^6$ & Iter & Error & Iter of Bica & Error of Bica \\ \hline
        10 & 1.0000e-06 & 16 & 5.2823e-06 & 17 & 3.964e-05 \\ 
        20 & 1.5625e-08 & 16 & 1.1847e-07 & ~ & ~ \\ 
        40 & 2.4414e-10 & 16 & 2.2197e-09 & ~ & ~ \\ 
        80 & 3.8147e-12 & 16 & 3.5748e-11 & ~ & ~ \\ 
        100 & 1.0000e-12 & 16 & 9.0258e-12 & 17 & 4.038e-009 \\ 
        200 & 1.5625e-14 & 17 & 7.4163e-14 & ~ & ~ \\  \hline
    \end{tabular}
\end{table}

\noindent {\bf Example 4.3.} (Example 4 in \cite{Bica1})
Consider the third order boundary value problem
\begin{align*}
u''' (t)&=e^{t}-\frac{1}{4}u(t)+\frac{1}{4}[u(\frac{t}{2}]^2, t\in [0, 1]\\
u(0)&=1,\ u'(0)=1, u'(1)=e
\end{align*}
with the exact solution $u(t)=e^{t}$. \\
The results of computation by Method 1 for Example 3 are presented in Table \ref{Tab3}.
For $n=10, 100, 1000$ the errors of approximate solution are 1.5160e-06, 1.5163e-10, 1.5099e-14, and the numbers of iterations are $8,8,9$, respectively.\\
\begin{table}[!ht]
    \centering
     \caption{The convergence of Method 1 for Example 4.3 }
     \label{Tab3}
    \begin{tabular}{ccccc}
        \hline
       $ n$ & $K$ & $Error$ & $h^4$ & $Order$ \\ \hline
        8 & 8 & 3.7008e-06 & 2.4414e-04 & ~ \\ 
        16 & 8 & 2.3135e-07 & 1.5259e-05 & 3.9997 \\ 
        32 & 9 & 1.4461e-08 & 9.5367e-07 & 3.9999 \\ 
        64 & 8 & 9.0379e-10 & 5.9605e-08 & 4.0000 \\ 
        128 & 9 & 5.6488e-11 & 3.7253e-09 & 4.0000 \\ 
        256 & 9 & 3.5305e-12 & 2.3283e-10 & 4.0000 \\ 
        512 & 9 & 2.2071e-13 & 1.4552e-11 & 3.9996 \\ 
        1024 & 9 & 1.4211e-14 & 9.0949e-13 & 3.9571 \\ 
    \end{tabular}
\end{table}

The convergence of  Method 2 for Example  4.3 is reported in Table \ref{Tab3a}, from which it is seen that the accuracy of the method is $O(h^6).$ 

\begin{table}[!ht]
    \centering
        \caption{The convergence of Method 2 for Example 4.3 }
     \label{Tab3a}
    \begin{tabular}{cccccc}
    \hline
 $N$ & $h^6$ & Iter & Error & Iter of Bica & Error of Bica \\ \hline
        10 & 1.0000e-06 & 5 & 6.9414e-09 & 8 & 1.4786e-06 \\ 
        20 & 1.5625e-08 & 5 & 6.8387e-11 & ~ & ~ \\ 
        40 & 2.4414e-10 & 4 & 6.3052e-12 & ~ & ~ \\ 
        80 & 3.8147e-12 & 4 & 4.3388e-13 & ~ & ~ \\ 
        100 & 1.0000e-12 & 4 & 1,7941e-13 & 8 & 1.4796e-10 \\ 
        200 & 1.5625e-14 & 4 & 1.1990e-14 & ~ & ~ \\ 
        300 & 1.3717e-15 & 3 & 8.8818e-16 \\ \hline
    \end{tabular}
\end{table}

\noindent {\bf Example 4.4.} (Example 4 in \cite{Bica})
Consider the third order boundary value problem
\begin{align*}
u''' (t)&=e^{-t}[u(t)]^{3/2}u(\frac{1}{2}t), t\in [0, 1]\\
u(0)&=1,\ u'(0)=1, u(1)=e.
\end{align*}
Notice that the boundary conditions in these examples are different from the ones of the previous examples. The Green function in this case is
\begin{equation}\label{eqG0a}
\begin{aligned}
G(t,s)=\left\{\begin{array}{ll}
-\frac{1}{2}(1-t)(2t-ts-s), \quad 0\le s \le t \le 1,\\
\, \, -\frac{1}{2}t^2(1-s)^2, \quad 0\le t \le s \le 1.\\
\end{array}\right.
\end{aligned}
\end{equation}
The results of computation for the example are given in Table \ref{TabExam5}.

\begin{table}[!ht]
    \centering
         \caption{The convergence of Method 2 for Example 4.4 }
     \label{TabExam5}
    \begin{tabular}{ccccc}
    \hline
        n & $h^6$ & iter & E & Order \\ \hline
        8 & 3.8147e-06 & 12 & 5.1727e-07 & ~ \\ 
        16 & 5.9605e-08 & 12 & 1.6924e-08 & 4.9338 \\ 
        64 & 1.4552e-11 & 12 & 1.7049e-11 & 4.9855 \\ 
        128 & 2.2737e-13 & 13 & 5.4001e-13 & 4.9806 \\ \hline
    \end{tabular}
\end{table}
From Table \ref{TabExam5} we see that the error is of $O(h^6)$ but the order (or rate) of convergence is only 5.
The comparison of the accuracy of Method 2 with \cite{Bica} is reported in Table 
\begin{table}[!ht]
    \centering
          \caption{The comparison of Method 2 in Example 4.4 with \cite{Bica} }
     \label{TabExam5a}
    \begin{tabular}{ccccc}
    \hline
        N & 10 & 100 & 800 & 1000 \\ \hline
        Method 2 & 1.7274e-07 & 1.8372e-12 & 6.6613e-16 & 4.4409e-16 \\ \hline
        Bica \cite{Bica} & 7.053e-07 & 6.703e-11 & ~ & 6.884e-15 \\ \hline
    \end{tabular}
\end{table}

\noindent {\bf Example 4.5.} (Example 2 in \cite{Bica})
Consider the third order boundary value problem
\begin{align*}
u''' (t)&=-\frac{4}{(1+t)^4}-( [u(t)]^4 +[u(t)]^3,). u(\frac{t}{2}) t\in [0, 1]\\
u(0)&=1,\ u(1)=\frac{1}{2}, u'(1)=-\frac{1}{4}
\end{align*}
with the exact solution $u(t)= \frac{1}{1+t}$.
The sixth order method gives the following results 

\begin{table}[!ht]
    \centering
          \caption{The convergence of Method 2 for Example 4.5 }
     \label{TabExam6}
    \begin{tabular}{ccccc}
    \hline
        N & $h^6$ & Iter  & E & Order \\ \hline
        8 & 3.8147e-06 & 13 & 2.0640e-06 & ~ \\ 
        16 & 5.9605e-08 & 13 & 5.1231e-08 & 5.3323 \\ 
        32 & 9.3132e-10 & 13 & 9.7636e-10 & 5.7135 \\ 
        64 & 1.4552e-11 & 14 & 1.6980e-11 & 5.8455 \\ \hline
    \end{tabular}
\end{table}

\begin{table}[!ht]
    \centering
             \caption{The comparison of Method 2 in Example 4.5 with \cite{Bica} }
     \label{TabExam6a}
    \begin{tabular}{ccccc}
    \hline
        N & 10 & 100 & 400 & 1000 \\ \hline
        Method 2  & 6.5312e-07 & 3.4605e-12  & 1.5432e-14  & 5.5511e-16  \\ 
        Bica \cite{Bica} & 2.1877e-06 & 2.402558e-10 & ~ & 2.398e-14 \\ \hline
    \end{tabular}
\end{table}

\noindent {\bf Example 4.6.} (Example 3 in \cite{Bica})
Consider the third order boundary value problem
\begin{align*}
u''' (t)&=-\frac{2}{3} u(t)-\frac{1}{3e^{-0.5t}} u(\frac{1}{2}), t\in [0, 1]\\
u(0)&=1,\ u(1)=\frac{1}{e}, u'(1)=-\frac{1}{e}
\end{align*}
with the exact solution $u(t)= e^{-t}$.
The sixth order method gives the following results 

\begin{table}[!ht]
    \centering
          \caption{The convergence of Method 2 for Example 4.6 }
     \label{TabExam7}
    \begin{tabular}{ccccc}
    \hline
        N & $h^6$ & Iter  & E & Order \\ \hline
        8 & 3.8147e-06 & 10 & 3.0241e-08 & ~ \\ \hline
        16 & 5.9605e-08 & 11 & 1.2083e-09 & 4.6454 \\ 
        32 & 9.3132e-10 & 11 & 4.1760e-11 & 4.8547 \\ 
        64 & 1.4552e-11 & 11 & 1.3666e-12 & 4.9334 \\ \hline
    \end{tabular}
\end{table}

\begin{table}[!ht]
    \centering
             \caption{The comparison of Method 2 in Example 4.6 with \cite{Bica} }
     \label{TabExam7a}
    \begin{tabular}{ccccc}
    \hline
        N & 10 & 100 & 500 & 1000 \\ \hline
      Method 2  & 1.1025e-08 & 1.7941e-13 & 3.3307e-16 & 2.2204e-16           \\ 
        Bica  & 6.7323e-08 & 6.5961e-12 & ~ & 7.772e-16 \\ \hline
    \end{tabular}
\end{table}

\section{Construction of fourth and sixth order accuracy for fourth order FDE}
In this section analogously as done for third order FDE in the previous section we construct methods of high order         accuracy for fourth order FDE equation. For simplicity we consider the problem on the interval [0, 1]:
\begin{equation}
\begin{split}
u^{(4)}(t)&= f(t,u(t),u(\varphi (t)), \; t\in (0, 1)\\
u(0)&=c_1; u(1)=c_2; u'(0)=c_3; u'(1)=c_4
\end{split}
\end{equation}
The Green function for this problem is
\begin{equation*}
	G(t,s)=\dfrac{1}{6}\begin{cases}
		s^2(t-1)^2(3t-s-2ts), \quad 0 \leq s \leq t \leq 1,\\
		t^2(s-1)^2(3s-t-2ts), \quad 0 \leq t \leq s \leq 1,
		\end{cases}
\end{equation*}
Discretizing the iterative method on the continuous level similar to  \eqref{iter1c}-\eqref{iter3c} we propose the following discrete iterative methods: Method 1 and Method 2 as follows;\\
{\bf Iterative Method 1 (Fourth order method): }\\
\begin{enumerate}[(i)]
\item Given 
\begin{equation}\label{iter1d4}
\Psi_0(t_i)=f(t_i,0,0),\ i=1,...,n+1. 
\end{equation}
\item Knowing $\Psi_k(t_i),\; k=0,1,...; \; i=1,...,n+1, $  compute approximately the definite integrals \eqref{iter2c} by the trapezoidal formulas with corrections:
\begin{equation}\label{iter2d4}
\begin{split}
U_k(t_i) &= g(t_i) + L_4(G,t_i)  \Psi_k, \\
V_k(t_i) &= g(\xi_i) + L_4(G,\xi_i)  \Psi_k,   i=1,...,n+1,
\end{split}
\end{equation}
where $\xi_i = \varphi(t_i)$ and $L_4(G,t_i) $  is defined on the grid function $\psi$ by the formula
\begin{equation*}
L_4(G, t_i) \psi = \sum_{j=1}^{n+1}h\rho_j G(t_i,t_j)\psi_j , \; i=1,...,n+1,
\end{equation*}
\item Update
\begin{equation}\label{iter3d4}
\Psi_{k+1}(t_i) = f(t_i,U_k(t_i),V_k(t_i)),\;  i=1,...,n+1.
\end{equation}
\end{enumerate}
{\bf Iterative Method 2 (Sixth order method): }\\
Let $m=[\frac{\xi_i}{h}]+1$. 
Method 2 is the same as Method 1 with a difference that in the step (ii) instead \eqref{iter2d4} there is
\begin{equation}\label{iter2d42}
\begin{split}
U_k(t_i) &= g(t_i) + L_6(G,t_i)  \Psi_k, \\
V_k(t_i) &= g(\xi_i) + L_6(G,\xi_i)  \Psi_k,   i=1,...,n+1,
\end{split}
\end{equation}
where for a grid function $\psi$
\begin{align*}
&L_6(G,t_i) \psi = 0, \; i=1, n+1\\
& 	L_6(G,t_i)\psi = \sum_{j=1}^{n+1}h\rho_j G(t_i,t_j)\psi_j +\dfrac{h^4}{720}\Big[ (3t_i^2-2t_i^3)\psi_{n+1}
 	+3t_i^2(1-t_i)D_2^{(1)}\psi_ {n+1}\\
 	&-(3t_i^2-2t_i^3-1)\psi_1
 	-3t_i(t_i-1)^2D_2^{(1)}\psi _1 - \psi_i\Big], \quad i=2,3,...,n.
\end{align*}
and
\begin{align*}
&L_6(G,\xi_i) \psi = 0, \; i=1\\
& 	L_6(G_0,\xi_i)\psi = \sum_{j=1}^{n+1}h\rho_j G(\xi_i,t_j)\psi_j +\dfrac{h^4}{720}\Big[ (3\xi_i^2-2\xi_i^3)\psi_{n+1}
 	+3\xi_i^2(1-\xi_i)D_2^{(1)}\psi_ {n+1}\\
 	&-(3\xi_i^2-2\xi_i^3-1)\psi_1
 	-3\xi_i(\xi_i-1)^2D_2^{(1)}\psi _1 - \psi_m -(\psi_{m+1}-\psi_m )(\xi_i-t_m)/h   \Big],\\
 	& \quad i=2,3,...,n+1.
\end{align*}
Analogously for the third order FDEs in the previous section, the described above Iterative methods 1 and 2 also have the accuracy $O(h^4)$ and $O(h^6)$. Below are some examples for demonstrating this assertion.
\section{Examples for  fourth order FDEs}
\noindent {\bf Example 5.1.} (Example 8, Problem (84) in \cite{Bica}) Consider the problem
\begin{equation}\label{eqExam2.O4}
	\begin{split}
		u^{(4)}(t)= e^{-t} [u(t)]^{3/2} u(\dfrac{1}{2}t),\; t\in (0,1),\\
		u(0)=1, \ u(1)= e, \ u'(0)=1, \ u'(1)=e
	\end{split}
\end{equation}
for which the exact solution is $u(t)=e^{t}$.\\
As was shown in \cite{Dang3} the problem has a unique solution. 
\begin{table}[!ht]
    \centering
    \caption{The convergence of Method 1 for \eqref{eqExam2.O4} }
     \label{Tab1.O4}
    \begin{tabular}{cccc}
    \hline
        n & K & E & Order \\ \hline
        8 & 9 & 4.5470e-07 & ~ \\ 
        16 & 9 & 2.8458e-08 & 3.9980 \\
        32 & 10 & 1.7940e-09 & 3.9876 \\ 
        64 & 9 & 1.1213e-10 & 3.9999 \\ 
        128 & 9 & 7.0077e-12 & 4.0001 \\ 
        256 & 9 & 4.3743e-13 & 4.0018 \\ \hline
    \end{tabular}
\end{table}
From Table \ref{Tab1.O4} it is obvious that the order of convergence of Method 1 is four.
\begin{table}[!ht]
    \centering
     \caption{The comparison of accuracy of Method 1 for \eqref{eqExam2.O4} with Bica method \cite{Bica} }
     \label{Tab2.O4}
    \begin{tabular}{cccc}
    \hline
        N & 10 & 100 & 1000 \\ \hline
        Method 1 & 1.8709e-07 & 1.8813e-11 & 1.7764e-15 \\ 
        Bica & 1.870948e-07 & 1.873168e-11 & 1.776357e-15 \\ \hline
    \end{tabular}
\end{table}
The convergence of Method 2 for \eqref{eqExam2.O4} is given in Table \ref{Tab2.O4}, from which it is seen that the accuracy of the method is $O(h^6)$.

\begin{table}[!ht]
    \centering
    \caption{The convergence of Method 2 for \eqref{eqExam2.O4} }
     \label{Tab2.O4}
    \begin{tabular}{ccccc}
    \hline
        N & $h^6$ & Iter  & E & Order \\ \hline
        8 & 3.8147e-06 & 7 & 1.7910e-09 & ~ \\ 
        16 & 5.9605e-08 & 7 & 2.8827e-11 & 5.9572 \\ 
        32 & 9.3132e-10 & 8 & 3.1861e-13 & 6.4995 \\ 
        64 & 1.4552e-11 & 8 & 2.5591e-14 & 3.6381 \\ 
        128 & 2.2737e-13 & 8 & 1.8874e-15 & 3.7612 \\ \hline
    \end{tabular}
\end{table}

\noindent {\bf Example 5.2.}
Now consider the problem (problem (85) in \cite{Bica})
\begin{equation}\label{eqExam2.O4a}
	\begin{split}
		u^{(4)}(t)=-4 e^{-t}+\frac{1}{2} u(t)+ e^{-t/2}u(\frac{1}{2}t),\; t\in (0,1),\\
		u(0)=1, \ u(1)= e^{-1}, \ u'(0)=1, \ u'(1)=0
	\end{split}
\end{equation}
for which the exact solution is $u(t)=te^{-t}$.\\
\begin{table}[!ht]
    \centering
    \caption{The convergence of Method 1 for \eqref{eqExam2.O4a} }
     \label{Tab1.O4a}
    \begin{tabular}{cccc}
    \hline
        N & K & E & Order \\ \hline
        8 & 7 & 8.5791e-07 & ~ \\  
        16 & 8 & 5.4197e-08 & 3.9845 \\ 
        32 & 8 & 3.4052e-09 & 3.9924 \\ 
        64 & 8 & 2.1285e-10 & 3.9998 \\ 
        128 & 8 & 1.3310e-11 & 3.9993 \\  \hline
    \end{tabular}
\end{table}
The comparison of accuracy of Method 1 for \eqref{eqExam2.O4a} with Bica method \cite{Bica} is given in Table \ref{Tab1.O4a}.
\begin{table}[!ht]
    \centering
     \caption{The comparison of accuracy of Method 1 for \eqref{eqExam2.O4a} with Bica method }
     \label{Tab2.O4a}
    \begin{tabular}{cccc}
    \hline
        N & 10 & 100 & 1000 \\ \hline
        Method 1 & 3.5625e-07 & 3.5727e-11 & 3.6082e-15 \\ 
        Bica & 3.224e-07 & 3.307e-11 & 3.331e-15 \\ \hline
    \end{tabular}
\end{table}

The convergence of Method 2 for \eqref{eqExam2.O4a} is given in Table \ref{Tab2.O4a}, from which it is seen that the accuracy of the method is $O(h^6)$.
\begin{table}[!ht]
    \centering
     \caption{The convergence of Method 2 for \eqref{eqExam2.O4} }
     \label{Tab3.O4a}
    \begin{tabular}{ccccc}
    \hline
        n & $h^6$ & iter & E & Order \\ \hline
        8 & 3.8147e-06 & 9 & 6.7167e-10 & ~ \\ 
        16 & 5.9605e-08 & 9 & 9.6914e-12 & 6.1149 \\ 
        32 & 9.3132e-10 & 9 & 2.6068e-13 & 5.2163 \\ 
        64 & 1.4552e-11 & 9 & 2.1094e-14 & 3.6274 \\ 
        128 & 2.2737e-13 & 9 & 1.9984e-15 & 3.3999 \\ \hline
    \end{tabular}
\end{table}

\section{On high order numerical methods for fifth order FDE}
In the previous sections we constructed the fourth and sixth order numerical methods for solving third and fourth order FDEs. In a similar way we can construct high order methods for fifth order FDEs. In this section we consider sixth order method for the following fifth order FDE
\begin{equation}\label{eqFDE5}
\begin{split}
&u^{(5)}(t) =f(x,u(t),u(\varphi (t))),\; t\in (0,1),\\
&u(0)=c_1 ,  u'(0) = c_2, \; u''(0)=c_3,\\
&u(1)= c_4, \; u'(1)=c_5.
\end{split}
\end{equation}
It is easy to verify that the Green function associated with the above problem is
\begin{equation}\label{GreenODE5}
\begin{aligned}
G(t,s)=\frac{1}{24}\left\{\begin{array}{ll}
s^2(t-1)^2(3s^2t^2+2s^2t+s^2-8st^2-4st+6t^2, &\quad 0\le s \le t \le 1,\\
(s-1)^3t^3(t-4s+3ts), &\quad 0\le t \le s \le 1.\\
\end{array}\right.
\end{aligned}
\end{equation}
We omit the formulas describing the sixth order numerical method for solving the  BVP \eqref{eqFDE5} and only report the results of  its convergence on two examples.\\
\noindent {\bf Example 7.1.} Consider the problem
\begin{equation}\label{Exam1FDE5} 
\begin{split}
&u^{(5)}(t) =\frac{1}{2}u(t) +\frac{1}{2}e^{\frac{3}{4}t}u(\frac{1}{4}t),\; t\in (0,1),\\
&u(0)=1 ,  u'(0) = 1, \; u''(0)=1,\\
&u(1)= e, \; u'(1)=e
\end{split}
\end{equation}
with the exact solution $u(t)= e^{t}$. The results of convergence of the sixth order method are given in Table \ref{Tab1.FDE5}:
\begin{table}[!ht]
    \centering
     \caption{The convergence of Method 2 for \eqref{Exam1FDE5} }
     \label{Tab1.FDE5}
    \begin{tabular}{ccccc}
    \hline
        n & $h^6$ & iter & E & Order \\ \hline
        8 & 3.8147e-06 & 6 & 4.5770e-10 & ~ \\ 
        16 & 5.9605e-08 & 6 & 7.1845e-12 & 5.9934 \\ 
        32 & 9.3132e-10 & 6 & 1.1169e-13 & 6.0073 \\ 
        64 & 1.4552e-11 & 6 & 1.7764e-15 & 5.9744 \\ 
    \end{tabular}
\end{table}

\noindent {\bf Example 7.2.} Consider the problem
\begin{equation}\label{Exam2FDE5} 
\begin{split}
&u^{(5)}(t) =720t+\frac{1}{5}u(t) u(t^2)-\frac{1}{5}(u(t))^3,\; t\in (0,1),\\
&u(0)=0 ,  u'(0) = 0, \; u''(0)=0,\\
&u(1)= 1, \; u'(1)=6
\end{split}
\end{equation}
with the exact solution $u(t)= t^6$. The results of convergence of the sixth order method are given in Table \ref{Tab2.FDE5}.
\begin{table}[!ht]
    \centering
     \caption{The convergence of Method 2 for \eqref{Exam2FDE5} }
     \label{Tab2.FDE5}
    \begin{tabular}{ccccc}
    \hline
        n & $h^6$ & iter & E & Order \\ \hline
        8 & 3.8147e-06 & 2 & 7.1262e-13 & ~ \\ 
        16 & 5.9605e-08 & 2 & 7.1748e-15 & 6.6341 \\ 
        32 & 9.3132e-10 & 2 & 1.2490e-16 & 5.8441 \\ 
        64 & 1.4552e-11 & 2 & 2.7756e-17 & 2.1699 \\ \hline
    \end{tabular}
\end{table}

By analogy with the fourth and sixth orders numerical methods, using the Euler-Maclaurain expansion \eqref{Euler-Mac} for $p=4$ we can construct an eighth order method named  Method 3 for the problem \eqref{eqFDE5}.
The numerical experiments on the above examples 7.1 and 7.2 show that indeed the method has the accuracy $O(h^8)$. See Tables \ref{Tab3.FDE5} and \ref{Tab4.FDE5}. From  these tables it is seen that the accuracy of the method is $O(h^8)$.
\begin{table}[!ht]
    \centering
     \caption{The convergence of Method 3 for \eqref{Exam1FDE5} }
     \label{Tab3.FDE5}
    \begin{tabular}{ccccc}
    \hline
        n & $h^8$ & iter & E & Order \\ \hline
        8 & 5.9605e-08 & 6 & 1.4213e-11 & ~ \\ 
        16 & 2.3283e-10 & 6 & 5.9064e-14 & 7.9107 \\ 
        32 & 9.0949e-13 & 6 & 8.8818e-16 & 6.0553 \\ 
        64 & 3.5527e-15 & 6 & 8.8818e-16 & 0 \\ \hline
    \end{tabular}
\end{table}

\begin{table}[!ht]
    \centering
     \caption{The convergence of Method 3 for \eqref{Exam2FDE5} }
     \label{Tab4.FDE5}
    \begin{tabular}{ccccc}
    \hline
        n & $h^8$ & iter & E & Order \\ \hline
        8 & 5.9605e-08 & 2 & 7.0173e-13 & ~ \\ 
        16 & 2.3283e-10 & 2 & 7.1887e-15 & 6.6090 \\ 
        32 & 9.0949e-13 & 2 & 1.2490e-16 & 5.8469 \\ 
        64 & 3.5527e-15 & 2 & 2.7756e-17 & 2.1699 \\ \hline
    \end{tabular}
\end{table}

\section{Concluding remarks}
In this paper we constructed high order numerical methods for solving the third, fourth and fifth orders FDE. From the results of computation for the examples it is possible to remark that  Method 2 (the sixth order method) applied to the fourth order FDE gives better accuracy than applied to the third order FDE and gives worse accuracy than applied to the fifth order FDE. This can be explained that the construction of the method requires  integration of a smooth function multiplied by the Green functions for the third, fourth and fifth orders equations, which have increasing smoothness. Thus, the higher the order of FDE, the more accurate the method can be constructed. This is also consistent with Bica's remark about the maximal  order of convergence of his method for solving functional differential equations because his method also leads to integration of the Green's function  multiplied by  the right-side function on each iteration. \par 
From the results of experiments for solving FDEs  we observe that the achieved  accuracy is somewhat worse than the accuracy  for ODEs in \cite{Dang4,Dang5,Dang6}. The reason of this may be in the computation of the values $\psi' (\xi_i) $, where $\xi_i =\varphi (t_i)$ in general not coinciding with grid points $t_i$ in Method 2.\par 
The technique used in this paper for the third, fourth and fifth orders FDE can be applied for designing numerical methods of order 8 and higher for FDE of sixth and higher orders.

\end{document}